\documentclass[12pt]{article}
\usepackage{amsmath, amsthm, amssymb}
\numberwithin{equation}{section}

\begin{document}
\author{Ajai Choudhry}
\title{An ancient diophantine problem of Brahmagupta}
\maketitle{}
\abstract{In the seventh century CE, Brahmagupta gave a parametric solution to the problem of finding two numbers $x$ and $y$ such that $x+y, x-y$ and $xy+1$ are all squares. We describe a method for generating infinitely many rational parametric solutions of the problem and another method for generating parametric solutions of the problem in terms of univariate polynomials with integer coefficients. The complete solution of the problem still remains elusive.} 



\section{Introduction}\label{intro}There are many interesting  diophantine problems dating back to ancient times for which the complete solution is still not known. This paper is concerned with one such problem considered by Brahmagupta in the seventh century CE and mentioned in the  partial English translation of his works by Colebrooke \cite[pp. 369-370, Art. 80-81]{Co}. The problem is to find  two numbers $x$ and $y$ such that $x+y, x-y$ and $xy+1$ are all squares. Brahmagupta gave the following rational solution of the problem:
\[
x=2a^2(a^2+b^2)/b^4, \quad y=2a^2(a^2-b^2)/b^4,
\]
where $a$ and $b$ are arbitrary parameters. As a numerical example, taking $(a, b)=(4, 2)$, he  obtained the solution $(x, y)=(40, 24)$. Taking $b=1$  yields infinitely many solutions of the problem in integers. With $b=1$ and assigning to $a$ the values  3 and 4, in succession, we get  two additional solutions with both $x$ and $y$ being positive integers not exceeding 1000 namely, $(x, y)=(180, 144) $ and $(x, y) = (544, 480)$.

While Dickson's {\it History} \cite[p.\ 444]{Di} mentions Brahmagupta's solution of this problem, there seems to be no further discussion of the problem in the mathematical literature. In this paper we obtain  infinitely many parametric solutions of the problem in rational numbers as well as in integers. 

\section{Parametric solutions of Brahma-\\gupta's  problem}\label{parmsols}
In Section \ref{rationalsols} we describe a method of generating infinitely many rational parametric solutions of Brahmagupta's problem using elliptic curves while in Section \ref{integersols} we use a generalized polynomial Pell's equation to  obtain infinitely many parametric solutions of the problem in terms of univariate polynomials with integer coefficients.
 
\subsection{Rational solutions}\label{rationalsols} 
If $x$ and $y$ are two rational numbers whose sum and difference are both squares, we may take $x=t^2+u^2, y=2tu$ so that $x+y=(t+u)^2$ and $x-y=(t-u)^2$. If $xy+1$ is also a square, there must exist a rational number $v$ such that the following condition is satisfied:
\begin{equation}
v^2 =2tu^3 + 2t^3u + 1. \label{condtuv}
\end{equation}
While Eq. \eqref{condtuv} is in three variables $t, u$, and $v$, and hence represents an algebraic surface over the field $\mathbb{Q}$, we may consider it as a cubic equation in the variables $u$ and $v$, and regard it as representing an elliptic curve over the function field $\mathbb{Q}(t)$. The  
birational transformtion defined by
\begin{align}
u &=X/(2t),  & v& =Y/(2t), \label{biratuv} \\
X & =2tu,  &Y &=2tv, \label{biratXY}
\end{align}
reduces  Eq. \eqref{condtuv} to the following Weierstrass model of an elliptic curve over the function field $\mathbb{Q}(t)$:
\begin{equation}
Y^2 = X^3 + 4t^4X + 4t^2. \label{ecWeier}
\end{equation}

It is readily seen that there is a rational point $P$ on the curve \eqref{ecWeier} whose coordinates are $(0, 2t)$. When $t=1$, the curve \eqref{ecWeier} reduces to the elliptic curve
\begin{equation}
Y^2 = X^3 + 4X + 4, \label{ecspl}
\end{equation}
and the coordinates of a  rational point $Q$ on the curve \eqref{ecspl} corresponding to the point $P$ on the curve \eqref{ecWeier} are $(0, 2)$. Direct computation using the group law on the elliptic curve \eqref{ecspl} shows that the points $2Q, 3Q$, and $4Q$ on the curve are given by $(1, -3)$, $(24, 118)$,  and $(-23/36, 235/216)$, respectively. Since the point $4Q$ has rational co-ordinates, it follows from the Nagell-Lutz theorem \cite[p. 56]{Si} on elliptic curves that the point $Q$ on the curve \eqref{ecspl} is not a point of finite order. In fact, the rank of the elliptic curve \eqref{ecspl}, as determined by the software SAGE, is 1. We can thus find infinitely many rational points on the curve \eqref{ecspl} using the group law.

Since in the special case when $t=1$, the point $Q$ on the curve \eqref{ecspl} corresponding to the point $P$ on the curve \eqref{ecWeier} is not a point of finite order, it follows that the point $P$ on the curve \eqref{ecWeier} cannot be a point of finite order. We can thus generate infinitely many rational points on the curve \eqref{ecWeier} using the group law. Each of these rational points on the curve \eqref{ecWeier} will yield a corresponding rational point on the curve \eqref{condtuv}, and using the relations $x=t^2+u^2, y=2tu$, we can obtain infinitely many pairs of rational numbers $x$ and $y$ in terms of the parameter $t$ such that $x+y, x-y$, and $xy+1$ are all  squares. 

We will now compute a couple of rational points on the elliptic curve and use them to obtain the desired pairs of rational numbers $x$ and $y$. The point $2P$ on the curve  \eqref{ecWeier} is given by $(X, Y) = (t^6, -t(t^8+2))$ corresponding to which we get the point $(u, v) =(t^5/2, -(t^8+2)/2)$ on the curve \eqref{condtuv} and hence we get the pair of rational numbers $x$ and $y$ given by 
\begin{equation}
(x, y) = ( t^2(t^8 + 4)/4,  t^6),  \label{ratsol1}
\end{equation} 
such that $x+y, x-y$ , and $xy+1$ are all  squares. Similarly, the point $3P$ on the curve \eqref{ecWeier} given by 
\[
(X, Y) =(8(t^8+2)/t^{10}, 2(3t^{16}+24t^8+32)/t^{15}),
\]
yields the following solution for $x$ and $y$:
\begin{equation}
(x, y) =( (t^{24} + 16t^{16} + 64t^8 + 64)/t^{22},  8(t^8 + 2)/t^{10}). \label{ratsol2}
\end{equation}
 Taking $t=1$ in each of the two solutions \eqref{ratsol1} and \eqref{ratsol2}, we get the following two rational solutions for $x$ and $y$: $(5/4, 1)$ and $(145,24)$.

\subsection{Integer Solutions} \label{integersols} 
Let $x$ and  $y$ be integers with the desired  properties. Since $xy+1$ is a square, we write $xy+1=(p+1)^2$, so that $xy=p(p+2)$ where $p$ is an integer. Since $x-y$ is a perfect square,  $p(p+2)$ must have more integer factors, hence we write $p=qr$, when we get $xy=qr(qr + 2)$. Accordingly, we now take
\begin{equation}
x=r(qr + 2), \quad y=q, \label{valxy}
\end{equation}
and now the following conditions must be satisfied:
\begin{align}
r(qr + 2) + q = s^2, \label{condsum} \\
r(qr + 2) - q = t^2, \label{conddiff}
\end{align}
where $q, r, s$, and $t$ are integers.
From Eqs. \eqref{condsum} and \eqref{conddiff}, we get, on subtraction, $2q=s^2-t^2$. We therefore take $s=t+2v$, when we get
\begin{equation}
q=2tv + 2v^2. \label{valq}
\end{equation}
Substituting the value of $q$ given by \eqref{valq} in Eq. \eqref{conddiff}, we get, after transposing all terms to one side,
\begin{equation}
t^2 - 2(r^2 - 1)vt - 2r^2v^2 + 2v^2 - 2r =0. \label{conddiff1}
\end{equation}
On writing
\begin{equation}
t=u+(r^2 - 1)v,  \label{valt}
\end{equation}
Eq. \eqref{conddiff1} reduces to
\begin{equation}
u^2 - (r^4 - 1)v^2 = 2r, \label{Pellgen}
\end{equation}
and now, on taking $r=2m^2$, we get the generalized Pell's equation
\begin{equation}
u^2 - (16m^8 - 1)v^2 = 4m^2, \label{Pellgenm}
\end{equation}
which has a solution $(u, v) =(8m^5, 2m)$. 

We also note that the related classical Pell's equation
\begin{equation}
\xi^2 - (16m^8 - 1)\eta^2 = 1, \label{Pellspl}, 
\end{equation}
has a solution
\begin{equation}
\xi =4m^4, \quad \eta=1. \label{Pellsplknownsol}
\end{equation}

 We will now generate a sequence of solutions of Eq. \eqref{Pellgenm} by using the   following well-known  identity of Brahmagupta \cite[pp. 14, 17 and 204]{We}:
\begin{equation}
(x_1^2-dy_1^2)(x_2^2-dy_2^2)=(x_1x_2+dy_1y_2)^2-d(x_1y_2+x_2y_1)^2. \label{ident1}
\end{equation}
If $(u_n, v_n)$ is any solution of Eq. \eqref{Pellgenm}, on taking  $d=(16m^8 - 1), x_1=u_n, y_1=v_n, x_2=4m^4, y_2=1, $ it follows from the identity \eqref{ident1} that a new solution $(u_{n+1}, v_{n+1})$ of Eq. \eqref{Pellgenm} is given by
\begin{equation}
\begin{aligned}
u_{n+1} & = x_1x_2+dy_1y_2 & =& \; 4m^4u_n + (16m^8 - 1)v_n,  \\
v_{n+1} & = x_1y_2+x_2y_1 \quad   & = & \; u_n +4m^4v_n.
\end{aligned}
\label{newsol}
\end{equation}
Since one solution of Eq. \eqref{Pellgenm} is given by $(u, v) =(8m^5, 2m)$, we take $(u_1, v_1)=(8m^5, 2m)$ and by repeatedly applying the recurrence relations \eqref{newsol}, we get an infinite sequence $(u_n, v_n), n=1, 2, ....$, of solutions of Eq. \eqref{Pellgenm} in terms of the parameter $m$. The next two solutions of this sequence are as follows:
\begin{equation}
\begin{aligned}
(u_2, v_2) & = (64m^9 - 2m,  16m^5), \\
(u_3, v_3) & = (512m^{13} - 24m^5,  128m^9 - 2m).
\end{aligned}
\label{solsPellgenm}
\end{equation}

Using the relations \eqref{valq} and \eqref{valt}, and noting that we have taken $r=2m^2$, it now follows from \eqref{valxy} that the values of ine integers $x$ and $y$ are given by
\begin{equation}
x = 4m^2(8m^6v^2 + 2m^2uv + 1), \quad  y = 8m^4v^2 + 2uv, \label{valxyf}
\end{equation}
where $(u, v)$ is a solution of the generalized Pell's equation \eqref{Pellgenm}. Now using the sequence $(u_n, v_n), n=1, 2, ....$, of solutions of Eq. \eqref{Pellgenm} and the relations \eqref{valxyf}, we can generate an infinite sequence of parametric solutions of our problem in terms of the parameter $m$. The first three solutions of this sequence are as follows:
\begin{equation}
\begin{aligned}
(x, y) & = (4m^2(64m^8 + 1),  64m^6), \\
(x, y) &= (4m^2(4096m^{16} - 64m^8 + 1), 64m^6(64m^8 - 1), \\
(x, y) & = (4m^2(262144m^{24} - 12288m^{16} + 128m^8 + 1), \\
& \quad \quad  128m^6(2048m^{16} - 96m^8 + 1)).
\end{aligned}
\label{solxy}
\end{equation}
The first of the above three solutions may also be obtained by taking $t=2m$ in the solution \eqref{ratsol1}.

As  numerical examples, taking $m=1$ in each of the three solutions given by \eqref{solxy}, we get three pairs of integers $(x, y)$ namely,
\[
(260, 64), \quad (16132,4032), \quad (999940, 249984), 
\]
which are such that, in each case, $x+y, x-y$ and $xy+1$ are all squares.
\section {Concluding remarks} 
While we can  generate infinitely many parametric solutions of Brahmagupta's problem both in rational numbers as well as in integers, we still do not have the complete solution of the problem. A simple computer search for integer solutions with both $x$ and $y$ not exceeding 1000 yielded six solutions. This naturally included the three integer solutions mentioned in the introduction and the solutions obtained in Section \ref{parmsols} but, in addition, we  found  a new solution $(793, 168)$ which cannot be generated either by Brahmagupta's original solution or by the parametric solutions obtained in this paper. Thus, the complete solution of the problem remains elusive.

\noindent Ajai Choudhry, 13/4 A Clay Square, Lucknow - 226001, India.

\noindent E-mail address: ajaic203@yahoo.com
 

\begin{thebibliography}{9}

\bibitem{Co} H. T. Colebrooke, {\it Algebra, with arithmetic and mensuration, from the Sanscrit of Brahmegupta and Bhascara, English translation}, John Murray, London, 1817. 

\bibitem{Di} L. E. Dickson, {\it History of the theory of numbers, Vol. 2}, Chelsea Publishing Company, USA, 1992

\bibitem{Si} J. H. Silverman and J. Tate, {\it Rational points on elliptic curves}, Springer-Verlag, New York, 1992. 

\bibitem{We} A. Weil, {\it Number theory, an approach through history from Hammurapi to Legendre}, Birkh\"{a}user  Birkhauser, Boston, 1983.
\end{thebibliography}
\end{document}